\documentclass[a4paper,12pt]{article}
\usepackage[utf8]{inputenc}
\usepackage[english,french]{babel}
   \usepackage{hyperref,varioref}
   \usepackage{amsmath}
   \usepackage{amsfonts}
   \usepackage{amssymb}
   \usepackage{mathrsfs}
   \usepackage{graphicx}
   \usepackage{textcomp}
   \usepackage{etex}
   \usepackage[abs]{overpic}
   \usepackage{float}
   \usepackage{cases}
   \usepackage{prettyref}
   \usepackage{empheq}
   \usepackage{relsize,exscale}
\newcommand{\R}{\mathbb{R}}

\newtheorem{thm}{Theora}[section]
\newtheorem{Theo}[thm]{Theorem}

\newtheorem{Lem}[thm]{Lemma}

\newenvironment{theorem*}[1]{\smallskip\noindent{\bf #1.}\rm}{\medskip}

\newenvironment{Proof}{\smallskip\noindent{\bf Proof.}\rm}
{\hfill $\Box$\medskip}
\newenvironment{proof}{\smallskip\noindent{\bf Proof}\rm}
{\hfill $\Box$\medskip}
\renewcommand\({\left(}
\renewcommand\){\right)}
\renewcommand\[{\left[}
\renewcommand\]{\right]}

\newcommand\E{{\rm E}}

\newcommand\la{\lambda}

\newcommand\rh{\rho}
\newcommand\si{\sigma}

\newcommand{\be}{\begin{equation}}
\newcommand{\ee}{\end{equation}}
\newcommand{\ba}{\begin{array}}
\newcommand{\ea}{\end{array}}
\newcommand{\bea}{\begin{eqnarray*}}
\newcommand{\eea}{\end{eqnarray*}}
\newcommand{\bean}{\begin{eqnarray}}
\newcommand{\eean}{\end{eqnarray}}
\newcommand\D{{Dim}}

\newcommand\se{\sigma}
\newcommand\ro{\rho}

\makeatletter \@addtoreset{equation}{section}

\makeatother
\begin{document}
\title{ On the Simplicity of Eigenvalues of Two Nonhomogeneous Euler-Bernoulli Beams Connected by a Point Mass}
\author{Jamel Ben Amara \thanks{Department of Mathematics, Faculty of Sciences of Tunis,
 University of Tunis El Manar, Mathematical Engineering Laboratory, Tunisia; e-mail: {\url{
jamel.benamara@fsb.rnu.tn}}.}~~~~~~~~~~Hedi Bouzidi
\thanks{Department of Mathematics, Faculty of Sciences of Tunis,
 University of Tunis El Manar, Mathematical Engineering Laboratory, Tunisia; e-mail:
{\url{hedi.bouzidi@fst.utm.tn}}.}}
\date{}
\maketitle {\bf Abstract:} In this paper we consider a linear system
modeling the vibrations of two nonhomogeneous Euler-Bernoulli beams
connected by a point mass. This system is generated by the following
equations\bea
&&\rho(x)y_{tt}(t,x)+(\sigma(x)y_{xx}(t,x))_{xx}-(q(x)y_{x}(t,x))_{x}=0,~t>0,~x \in (-1,0)\cup(0,1),\\
&&M
y_{tt}(t,0)=\({T}y(t,x)\)_{\mid_{x=0^-}}-\({T}y(t,x)\)_{\mid_{x=0^+}},~~t>0,\eea
with hinged boundary conditions at both ends, where ${T}y =
(\se(x)y_{xx})_{x} - q(x)y_x$. We prove that all the associated
eigenvalues $\(\la_n\)_{n\geq1}$ are algebraically simple,
furthermore the corresponding eigenfunctions $\(\phi_n\)_{n\geq1}$
satisfy $\phi_n'{T}{\phi}_{n}(-1)>0$ and $\phi_n'{T}{\phi}_{n}(1)<0$
for all $n\geq1$. 
These results give a key to the solutions of
various control and stability problems related to this system.\\


\section{Introduction}
In the last three decades there has been an increasing interest in
the study of the dynamics and control of various hybrid models for
systems of rods, strings and beams with attached masses.  For the
boundary controllability and stability problems related to this type
of systems we can refer to \cite{Avd, JA, Beam2, Guo2001, ren2}, see
also \cite{
 Conrad1998, Conrad1994, E.Z} and references therein. As is well known, the spectral analysis is the key tools for solving
these problems.\\
\indent In this paper, we consider a one-dimensional linear hybrid
system which is composed by two nonhomogeneous hinged
Euler-Bernoulli beams connected by a point mass. We assume that the
first beam occupies the interval $\Omega_{1}=(-1, 0)$ and the second
one occupies the interval $\Omega_{2}=(0, 1)$. The vibrations of the
first and the second beam will be, respectively, presented by the
functions \bea y:=y(t,x),&(t,x) \in
\(0,\infty\)\times\Omega_{1}\hbox{ and } z:=z(t,x),&(t,x) \in
\(0,\infty\)\times\Omega_{2}. \eea The position of the mass $M> 0$
attached to the beams at the point $x = 0$ is denoted by the
function $y:= y(t,0)$ for $t > 0$. The dynamic behavior of the
system is governed by the following PDE: \bean\label{beam1} \left\{
\begin{array}{ll}
\rho_{1}(x)y_{tt}(t,x)+(\sigma_{1}(x)y_{xx}(t,x))_{xx}-(q_{1}(x)y_{x}(t,x))_{x}=0,&(t,x) \in \(0,\infty\)\times\Omega_{1},\\
\rho_{2}(x)z_{tt}(t,x)+(\sigma_{2}(x)z_{xx}(t,x))_{xx}-(q_{2}(x)z_{x}(t,x))_{x}=0,&(t,x) \in \(0,\infty\)\times\Omega_{2},\\
y(t,-1)=y_{xx}(t,-1)=z(t,1)=z_{xx}(t,1)=0,&t\in\(0,\infty\),\\
{y}(t,0)= {z}(t,0),{y_x}(t, 0)={z_x}(t, 0),\sigma_1y_{xx}(t, 0) =\sigma_2z_{xx}(t, 0), &t\in\(0,\infty\),\\
My_{tt}(t, 0)={{T}^1 y}(t, 0)-{{T}^2 z}(t, 0),&t\in\(0,\infty\),
\end{array}
\right. \eean where ${T}^if(t,x):= (\se_i(x)f_{xx}(t,x))_{x} -
q_i(x)f_x(t,x)$, for $t>0$ and $x\in \Omega_{i}~(i=1,2)$. The
coefficients $\rho_{i}$, $\sigma_{i}$ and $q_{i}~(i=1,2)$ of each
beam represent, the density, the flexural rigidity and the axial
force, respectively, see for instance \cite{Virgin, Roseau}. By
applying separation of variables to System \eqref{beam1}, we obtain
the following spectral problem :
\begin{empheq}[left=\empheqlbrace]{align}
&(\se_1(x)u'')''-(q_1(x)u')'=\la\ro_1(x)u,~~x \in \Omega_{1},\label{s1}\\
&(\se_2(x)v'')''-(q_2(x)v')'=\la\ro_2(x)v,~~x\in \Omega_{2},\label{s2}\\
&u(-1)=u''(-1)=v(1)=v''(1)=0,\label{s3}\\
&u(0)=v(0),~~u'(0)=v'(0),~~\se_1u''(0)=\se_2v''(0),\label{s4}\\
&\(\mathcal{T}^1u(x)-\mathcal{T}^2v(x)\)_{\mid_{x=0}}=-M\la
u(0),\label{s5}\end{empheq} where $\mathcal{T}^if(x):=
(\se_i(x)f''(x))' - q_i(x)f'(x)$ for $x\in \Omega_{i}~(i=1,2)$.
Throughout this paper we assume that \be \rho_i \in
C(\Omega_{i}),~\sigma_i \in H^{2}(\Omega_{i}),~q_i\in
H^{1}(\Omega_{i}), \label{beam2}\ee and there exist constants
$\rho_0,~\sigma_0>0$, such that \be
\rho_i(x)\geq\rho_{0},~~\se_i(x)\geq\se_{0},~~ q_i(x)\geq0,~~x \in
 \Omega_{i}~(i=1,2). \label{beam3}\ee
\indent There exists an extensive mathematical and engineering
literature devoted to the spectral analysis for various systems of
vibrating beams. The asymptotics, the simplicity of eigenvalues and
the oscillations of the eigenfunctions with their derivatives of
vibrating beams without point mass have been investigated in
\cite{BK, BK1, B22, B11, Janc} for different boundary conditions.
These results were extended in a number of works to Euler-Bernoulli
beams with end masses, see \cite{Aliev1, Aliev2, B1, B2, Kerimov}.
However, the spectral proprieties related to a series of beams with
interior attached masses have been considered only in the case of
constant physical parameters. Namely, by using a precise spectral
analysis together with the theory of non-harmonic Fourier series
Castro and Zuazua \cite{Beam1, Beam2, Beam3} proved the exact
controllability for two type of homogenous flexible beams connected
by a point masse. Later on, Mercier and R\'{e}gnier
\cite{ren1,ren2}, extended their results to the case of network of
Euler-Bernoulli beams with interior attached
masses.\\
\indent The main result of this paper is the following.
\begin{Theo}\label{th1}~~~~~~~~~~~~
\begin{description}
\item[(a)]The eigenvalues $\(\la_{n}\)_{n\geq1}$ of the spectral problem
\eqref{s1}-\eqref{s5} are real, algebraically simple and form an
infinitely increasing sequence such that
$$0<\lambda_{1}<\lambda_{2}<.......<\lambda_{n}<.....\underset{n\rightarrow \infty}{\longrightarrow}\infty .$$
 \item [(b)] The corresponding eigenfunctions
$({\phi}_{n})_{n\in\mathbb{N^{*}}}$ have the following properties:
\be {\phi}_{n}'\mathcal{T}^1{\phi}_{n}(-1)<0 \hbox{ and
}{\phi}_{n}'\mathcal{T}^2{\phi}_{n}(1)>0 ~~\hbox{for all
}n\geq1.\label{alt12} \ee
\end{description}
\end{Theo}
The proof of this Theorem is mainely based on some properties of
fourth-order linear differential equations (see \cite{LN}) and the
associated theory of subwronskians (see \cite{Barret1,Barret2}).\\
\indent This paper is organized as follows. In Section $2$, we
associate to Problem \eqref{s1}-\eqref{s5} a self-adjoint operator
with compact resolvent defined in a well chosen Hilbert space. In
Section $3$, we establish several lemmas that are used in the proof
of Theorem \ref{th1}, and that we believe are of independent
interest. Finally in Section $4$, we give the proof of Theorem
\ref{th1}.
\section{Operator framework}
Let us define the Hilbert space
$$
\mathcal{H}=L^{2}(-1,0)\times L^{2}(0,1)\times \mathbb{R},
$$
with the scalar product $\langle./.\rangle_{\mathcal{H}}$ defined
by: for all $y_{i}=(u_{i},v_{i},z_{i})^{t}\in \mathcal{H}~(i=1,2)$,
where $^{t}$ denotes the transposition, we have
$$
\langle y_{1},
y_{2}\rangle_{\mathcal{H}}=\int_{\Omega_{1}}u_{1}u_{2}\rho_{1}(x)dx+
\int_{\Omega_{2}}v_{1}v_{2}\rho_{2}(x)dx+Mz_{1}z_{2}.
$$
Let \bea \mathcal{V}=\{(u,v)\in H^{2}(\Omega_{1})\times
H^{2}(\Omega_{2})~ :\hbox{ satisfying } \eqref{s3}, \eqref{s4}\},
 \eea endowed with the norm $$
\|(u,v)\|^{2}_{\mathcal{V}}=\int_{\Omega_{1}}|u''(x)|^{2}dx+
\int_{\Omega_{2}}|v''(x)|^{2}dx.$$ Here and it what follows
$H^k(\Omega_{i})~(i=1,2),$ refers to the standard Sobolev space,
with $H^0 = L^2$. It is easy to see that $\mathcal{V}$ is
algebraically and topologically equivalent to $H^{2}\cap
H^{1}_{0}(-1,1)$. Let us consider the following closed subspace of
$\mathcal{V} \times \mathbb{R}$
$$
\mathcal{W}=\{(u,v,z)\in \mathcal{V} \times \mathbb{R}:
u(0)=v(0)=z\},
$$
equipped with the norm
$\|(u,v,z)\|^{2}_{\mathcal{W}}=\|(u,v)\|^{2}_{\mathcal{V}}$. We
introduce the operator $\mathcal{A}$ defined in $\mathcal{H}$ by
setting:
\begin{equation}\label{est}
\mathcal{A}y=\begin{cases} \frac{1}{\rho_{1}(x)}((\sigma_{1}(x)u'')''-(q_{1}(x)u')',\\
\frac{1}{\rho_{2}(x)}((\sigma_{2}(x)v'')''-(q_{2}(x)v')',\\
-\frac{1}{M}\(\mathcal{T}^1u(x)-\mathcal{T}^2v(x)\)_{\mid_{x=0}},
\end{cases}
\end{equation}
where $y=(u,v,z)^{t}$ on the domain \bea \mathcal D(\mathcal{A})=\{
                    (u,v,z)\in \mathcal{W}~:~(u,v)\in
H^{4}(\Omega_{1})\times H^{4}(\Omega_{2}) \} \eea which is dense in
$\mathcal{H}$. Obviously, Problem \eqref{s1}-\eqref{s5} is
equivalent to the following spectral problem
$$\mathcal{A}\phi=\la \phi,~~\phi=(u,v,z)^{t}\in D(\mathcal{A}),$$
i.e., the eigenvalues $\(\la_n\)_{n\geq1}$, of the operator
$\mathcal{A}$ and Problem \eqref{s1}-\eqref{s5} coincide together
with their multiplicities. Moreover, there is a one-to-one
correspondence between the eigenfunctions,
$$\phi_n(x)=(u_n(x),v_n(x),z_n)^{t}\leftrightarrow(u_n(x),v_n(x))^{t},~~z_n=u_n(0),~n\geq1.$$
\begin{Theo}\label{rr}
The linear operator $\mathcal{A}$ is positive and self-adjoint such
that $\mathcal{A}^{-1}$ is compact. Moreover, the spectrum of
$\mathcal{A}$ is discrete and consists of a
 sequence of positive eigenvalues
$(\lambda_{n})_{n\in\mathbb{N}^{*}}$ tending to $+\infty$:
$$0<\lambda_{1}\leq\lambda_{2}\leq.......\leq\lambda_{n}\leq.....
\underset{n\rightarrow +\infty}{\longrightarrow}+\infty.$$
\end{Theo}
\begin{Proof}
Let $y=(u_1,u_2,z)\in\mathcal D(\mathcal{A})$, then by integration
by parts, we have \bea \langle \mathcal{A}y,
y\rangle_{\mathcal{H}}&=&\sum_{i=1}^2\int_{\Omega_i}\Big{(}(\sigma_i(x)u_i'')''-
(q_i(x)u_i')'\Big{)}u_i dx-\(\mathcal{T}^1u_1(x)-\mathcal{T}^2u_2(x)\)_{\mid_{x=0}}z, \\
&=&\int_{\Omega_1}\(\sigma_1(x)|u_1''|^{2}+q(x)|u_1'|^{2}\)dx+
\int_{\Omega_2}\(\sigma_2(x)|u_2''|^{2}+q(x)|u_2'|^{2}\)dx \eea
Since $\sigma_i>0$ and $q_i\geq0 (i=1,2)$ then $\langle
\mathcal{A}y, y\rangle_{\mathcal{H}}>0,$ and hence, the linear
operator $\mathcal{A}$ is positive. Furthermore, it is easy to show
that $Ran(\mathcal{A}-iId)=\mathcal{H}$, and this implies that
$\mathcal{A}$ is selfadjoint. Since the space $\mathcal{W}$ is
continuously and compactly embedded in the space $\mathcal{H}$, then
$\mathcal{A}^{-1}$ is compact in $\mathcal{H}$. The proof is
complete.
\end{Proof}
\section{Basic Lemmas}
In this section, we establish several basic results that will be
used frequently in the next section. We consider the linear fourth
order differential equation defined on the interval $\[a,b\]$,
$a\geq0$: \be\label{ss*} (\sigma(x) u'')''-(q(x) u')'-\rho(x) u=0,
\ee where the functions $\rho(x)$, $\sigma(x)$ are uniformly
positive, and $q(x)$ is nonnegative such that
$$\rho \in C(a,b),~\sigma\in H^{2}(a,b),~q\in H^{1}(a,b).$$
We start by mentioning the following lemma due Leighton-Nehari
\cite{LN}.
\begin{Lem}\cite[Lemma 2.1]{LN}\label{Ls1}
Let $u$ be a nontrivial solution of the differential equation
\eqref{ss*} for $q\equiv0$.
 If $u, u',
u''$ and $(\sigma u '')'$ are nonnegative at $x=a$ (but not all
zero) they are positive for all $x>a$. If $u, -u', u''$ and
$-(\sigma u'')'$ are nonnegative at $x=b$ (but not all zero) they
are positive for all $x<b$.
\end{Lem}
The following lemma was stated in \cite[Lemma 2.1]{BK1}.  For the
reader's convenience, we propose here a simpler proof.
\begin{Lem}\label{Ls1m}
Let $u$ be a nontrivial solution of Equation \eqref{ss*}.
 If $u, u',
u''$ and $Tu= (\se(x)u'')' - q(x)u'$ are nonnegative at $x=a$ (but
not all zero), then they are positive for all $x>a$. If $u, -u',
u''$ and $\(-Tu\)$ are nonnegative at $x=b$ (but not all zero), then
they are positive for all $x<b$.
\end{Lem}
\begin{Proof}
Let $h$ be the unique solution of the following second order initial
value problem: \begin{empheq}[left=\empheqlbrace]{align}
&(\sigma(x) h')'- q(x)h=0,~~x\in (a, b]\label{sec}\\
&h(a)=1,~~h'(a)=0.\label{s11}
\end{empheq}
It is known, by Sturm comparison theorem \cite[Chapter 1]{BI} that
$h(x)>0$ on $[a,b]$. Hence, the following modified Leighton-Nehari
substitution \cite[Theorem 12.1]{LN}
$$
t(x):=\gamma^{-1}{(b-a)\int_a^x h(s)ds}+a,~~\gamma={\int_a^b
h(s)ds},
$$
transforms Equation \eqref{ss*} into
\begin{equation}\label{s1t}
{\(\tilde\si(t)\ddot{\tilde u}\)}^{..} =\tilde \rh(t) \tilde u,~~t\in[a,b],\\
\end{equation}
where $\tilde \si(t)=(\gamma (b-a)^{-1}
h(x(t)))^3\si(x(t)),~\tilde\rh(t)=\gamma (b-a)^{-1}
h^{-1}(x(t))\rh(x(t)) ~\hbox{ and }~  ^\cdot:=\frac{d}{dt}.$ If $u$
is a nontrivial solution of \eqref{ss*}, then $\tilde u(t)\equiv
u(x(t))$ is a nontrivial solution of \eqref{s1t}. Furthermore, we
have
\begin{equation}\label{s2t}
{\dot {\tilde u}}=\gamma (b-a)^{-1}h^{-1}u',~
 \gamma^2 (b-a)^{-2}h^{3}{\ddot {\tilde u}}=h u''- u'h',~
{\(\tilde\si{\ddot {\tilde u}}\)}^{.}=Tu.
\end{equation} It is easy to see from \eqref{s11} and \eqref{s2t},
that $u, u', u''$ and $Tu$ are positive at $x=a$. Hence, in view of
Lemma \ref{Ls1}, \be\tilde{u}>0,
~\dot{\tilde{u}}>0,~\ddot{\tilde{u}}>0,~{( \tilde\si\ddot
{\tilde{u}})}^{.}>0 \hbox{ in } (a,b].\label{s12}\ee Since $\sigma
h'(x)=\int_a^xqh\rh(x)dx$, then $h'(x)>0$ on $(a,b].$ Therefore,
combining \eqref{s2t} and \eqref{s12}, one gets
$$u>0, ~u'>0,~u''>0,~Tu
>0 \hbox{ in }(a,b].$$ For
the proof of the second statement it is sufficient to replace the
initial conditions \eqref{s11} by
$$h(b)=1,~~h'(b)=0.$$
By Sturm comparison Theorem, $h>0$ on the interval $[a,b]$. The
Lemma is proved.
\end{Proof}\\
Using Lemma \ref{Ls1m}, we can establish the following lemma.
\begin{Lem}\label{le:dim}~~ \begin{enumerate}
                            \item Let $\E_{u}$ be the space
                            of solutions of Equation \eqref{s1} for $\la>0$, satisfying one of the following sets of boundary
                            conditions :
\begin{empheq}[left=\empheqlbrace]{align}
&u(-1) = u''(-1)=0, ~~\alpha_1 u'(0)=\beta_1
u''(0),\label{s13}\\
&u(-1) = u''(-1)=0,~~ \alpha_1 \mathcal{T}^1u(0)=\beta_1
u(0),\label{s14}
\end{empheq}
where $(\alpha_1,\beta_1)\in \R^2\backslash\{(0,0)\}$ and
$\alpha_1\beta_1\leq0$. Then $\D\E_{u}=1$.
\item Let $\E_{v}$ be the space of solutions of Equation \eqref{s2} for $\la>0$, satisfying
one of the following sets of boundary conditions:
\begin{empheq}[left=\empheqlbrace]{align}
&v(1) = v''(1)=0, ~~\alpha_2 v'(0)=\beta_2
v''(0),\label{sv13}\\
& v(1) = v''(1)=0,~~ \alpha_2 \mathcal{T}^2v(0)=\beta_2
v(0),\label{sv14}
\end{empheq}
where $(\alpha_2,\beta_2)\in \R^2\backslash\{(0,0)\}$ and
$\alpha_2\beta_2\geq0$. Then $\D\E_{v}=1$.
                          \end{enumerate}
\end{Lem}
\begin{Proof} Suppose that there exist two linearly independent
solutions $u_i(i=1,2)$ of Problem \eqref{s1}-\eqref{s13} . Both $
u'_{1}(-1)$ and $ u'_{2}(-1)$ must be different from zero since
otherwise it would follow from Lemma~\ref{Ls1m} that $
u'_iu''_{i}(0)>0(i=1,2)$ and this is in contradiction with the last
boundary condition in \eqref{s13}. In view of the assumptions about
$u_1$ and $u_2$, the solution
 $$u(x)=  u'_{2}(-1)u_{1}(x) -  u'_{1}(-1)u_{2}(x)$$
  satisfies
 $u(-1)= u'(-1)= u''(-1)=0$ and $u'u''(0)\leq0$. This again contradicts
Lemma~\ref{Ls1m} unless $u\equiv0$. The other statements of the
Lemma can be proved in a same way.
\end{Proof}\\
\begin{Lem}\label{le:cont}
Every solution $\phi_0(x)$ of the regular problem
\eqref{s1}-\eqref{s5} for $M=0$ has only simple zeros in $(-1,1)$.
\end{Lem}
\begin{Proof}
Without loss of generality, assume that there exists
$x_{0}\in\Omega_1$ such that $\phi_0(x_0)=\phi'_0(x_0)=0$. If
$\phi''_0\mathcal{T}^1\phi_0(x_0)<0$, then the second statement of
Lemma \ref{Ls1m} yields a contradiction with the boundary condition
$\phi_0(-1)=0$. Now, if $\phi''_0\mathcal{T}^1\phi_0(x_0)\geq0$,
then $\phi''_0\mathcal{T}^1\phi_0(0)\geq0$, and hence from
\eqref{s4} with $M=0$, we have
$\phi''_0\mathcal{T}^2\phi_0(0)\geq0$. Therefore by the first
statement of Lemma \ref{Ls1m}, $\phi_0(1)\neq0$, a contradiction.
The proof is complete.
\end{Proof}\\
It is known that any solution of Equation \eqref{s1} which satisfies
the initial conditions $u(-1) = u''(-1)=0$ may be expressed as a
linear combination of ${y}_1(x)$ and ${y}_2(x)$, where
${y}_i~(i=1,2),$ are the fundamental solutions of \eqref{s1}
satisfying the initial conditions: \bean  &{y}_1(-1) = y''_1(-1) =
\mathcal{T}^1y_1(-1)=0,\quad y'_1(-1)=1,\label{eq:initu}\\ & y_2(-1)
= y'_2(-1)= y''_2(-1)=0,\quad \mathcal{T}^1 y_2(-1)=1.
\label{eq:initv}\eean In view of Lemma \ref{Ls1m}, ${y}_{i}$,
${y}'_{i}$, ${y}''_{i}$ and $\mathcal{T}^1 {y}_{i}$ $(i=1,2)$ are
positive in $\Omega_1\cup\{0\}$. We introduce the following
subwronskians (see \cite{Barret1, Barret2}):
\begin{equation}\label{eq:wrons1}
\overline\sigma_1= y_1y_2' -y_2y_1', \quad \overline\sigma'_1= y_1
y_2'' -y_2y_1'',\quad \overline\tau_1= y_1\mathcal{T}^1y_2 -
y_2\mathcal{T}^1y_1.
\end{equation}
Clearly, if for some $\lambda>0$ and $x_{0}\in\Omega_1\cup\{0\}$,
$\overline{\sigma}_1(x_{0}, \lambda)=0$, then $\lambda$ is an
eigenvalue and $u(x)=y_{1}(x_{0})y_{2}(x)-y_{2}(x_{0})y_{1}(x)$ is
the corresponding eigenfunction of the problem determined by
\eqref{s1} and the boundary conditions
\begin{equation}\label{s'1}
u(-1)=u''(-1)=u(x_{0})=u'(x_{0})=0.
\end{equation}
Similar conclusions for the other subwronskians $\overline{\tau}_1$
and $\overline{\sigma}_1'$.\\
Analogously, we introduce the following subwronskians associated
with Equation \eqref{s2} :
\begin{equation}\label{eq:wronss1}
\overline\sigma_2= z_1z_2' -z_2z_1', \quad \overline\sigma'_2=
z_1z_2'' -z_2z_1'',\quad \overline\tau_2=z_1\mathcal{T}^2z_2 -
z_2\mathcal{T}^2z_1 ,
\end{equation}
 where $z_1$ and $z_2$ are two linearly independent solutions of
\eqref{s1} which satisfy the initial conditions: \bean
&z_1(1)=z''_1(1)=\mathcal{T}^2 z_1(1)=0,~z'_1(1)=-1,\label{c1}\\
&z_2(1)=z'_2(1)= z''_2(1)=0,~\mathcal{T}^2z_2(1)=-1.\label{c2} \eean
Obviously, we have
\begin{equation}\label{c3}
z_i>0, ~z''_i>0, ~z'_i<0,~ \mathcal{T}^2z_i<0,~~~i=1,2.
\end{equation}
If one of the subwronskians $\overline{\sigma}_2$,
$\overline{\sigma}_2'$ and $\overline{\tau}_2$ vanishes for some
$\lambda>0$ and $x_{0}\in\Omega_2\cup\{0\}$, then $\lambda$ is an
eigenvalue of the problem determined by \eqref{s2} and the boundary
conditions
\begin{equation}\label{s'+++1}
v(x_{0})=v'(x_{0})=v(1)=v''(1)=0.
\end{equation}
\begin{Lem}The following formulas hold: \be\label{eq:wrons2}
\overline{\tau}_1={y}_{1}'\sigma_{1}{y}_{2}''-{y}_{2}'\sigma_{1}{y}_{1}''~~\hbox{
and
}~~\overline{\tau}_2={z}_{1}'\sigma_{2}{z}_{2}''-{z}_{2}'\sigma_{2}{z}_{1}''.
\ee
\end{Lem}
\begin{Proof} By multiplying Equation \eqref{s1} (for $u=y_1$) by
$y_2$, and twice integrating by parts from $-1$ to $x$, yields
$$
   y_1\mathcal{T}^1y_2(x)-\sigma_{1}y_1''y_2'(x)=
   \la \int_{-1}^{x} y_1y_2\rho_1(x)dx-
   \int_{-1}^{x} \(\sigma_{1}y_1''y_2''(x)+q_{1}y_1'y_2'(x)\)dx.
     $$
Similarly, $$
   y_2\mathcal{T}^1y_1(x)-\sigma_{1}y_2''y_1'(x)=
   \la \int_{-1}^{x} y_1y_2\rho_1(x)dx-
   \int_{-1}^{x} \(\sigma_{1}y_1''y_2''(x)+q_{1}y_1'y_2'(x)\)dx.
     $$
Subtracting these equalities together with \eqref{eq:wrons1}, we get
the expression \eqref{eq:wrons2}. The proof of the second expression
in \eqref{eq:wrons2} is similar.
\end{Proof}
\begin{Lem}\label{lem} Let $\la>0$ and fixed $x_{0}\in\Omega_i\cup\{0\}~(i=1,2)$. If one of the subwronskians
 $\overline{\sigma}_i$,
$\overline{\sigma}_i'$, and $\overline{\tau}_i~(i=1,2)$ vanishes at
$(x_{0}, \lambda)$, then all the other subwronskians are different
from zero at $(x_{0}, \lambda)$.
\end{Lem}
\begin{Proof} Without loss of generality, suppose that
$\overline{\sigma}_1(x_{0},\lambda)=\overline{\sigma}_1'(x_{0},\lambda)=0$,
for some $x_{0}\in\Omega_1\cup\{0\}$ and $\la>0$. Then there exist
two eigenfunctions $u_1$ and $u_2$ of the problems determined by
Equation \eqref{s1} and the boundary conditions \eqref{s'1} and
\begin{equation}\label{s'2}
u(-1)=u''(-1)=u(x_{0})=u''(x_{0})=0,
\end{equation}
respectively. Since $\varphi_{1}$ and $\varphi_{2}$ satisfy in
common $u(-1)=u''(-1)=u(x_{0})=0$, then by Lemma \ref{le:dim}, they
are colinear, i.e., $\varphi_1=c\varphi_2$, $c\in\R$. This means
that
$\varphi_{1}(x_{0})=\varphi'_{1}(x_{0})=\varphi''_{1}(x_{0})=0$, and
this is in contradiction with the second statement of Lemma
\ref{Ls1m}. The proof is complete.
\end{Proof}
\section{Proof of Theorem \ref{th1}}
In this section we prove Theorem \ref{th1}.\\
\begin{proof} {\bf of Assertion (a).}
Let $\la$ be an eigenvalue of Problem \eqref{s1}-\eqref{s5} and let
$E_{\la}$ be the corresponding eigenspace (i.e.,
$E_{\la}=Ker\(\mathcal{A}-\la I\)$). We shall prove that $\D
E_{\la}=1$. To this end, let $\varphi_1,~\varphi_2$ and $\varphi_3$
be three solutions of Equation \eqref{s1}, satisfying the boundary
conditions
 \bean &\varphi_1(-1)=\varphi_1''(-1)=\varphi_1(0)=0,\label{c4}\\
  &\varphi_2(-1)=\varphi_2''(-1)=\varphi_2'(0)=0,\label{c5}\eean
 and
 \be
\varphi_3(-1)=\varphi''_3(-1)=\varphi''_3(0)=0,\label{u}\\
 \ee
 respectively. By virtue of Lemma \ref{le:dim}, $\varphi_1,~\varphi_2$ and
$\varphi_3$, are the unique solutions, up to a multiplicative
constant, of Equation \eqref{s1} satisfying
 the boundary conditions \eqref{c4}, \eqref{c5} and \eqref{u}, respectively.
Similarly, let $\psi_i ~(i=1,2,3)$ be the unique solutions of
Equation
 \eqref{s2}, satisfying the boundary conditions
 \bean &\psi_1(0)=\psi_1(1)=\psi_1''(1)=0,\label{c6}\\
 &\psi_2'(0)=\psi_2(1)=\psi_2''(1)=0,\label{c7}\eean
 and \be \psi''_3(0)=\psi_3(1)=\psi''_3(1)=0,\label{v}\ee
 respectively. It is easy to see that $\varphi_i$ and $\psi_i~(i=1,2,3),$ can be written as follows:
\bean &&\varphi_1(x)= y_1(0)y_2(x)- y_2(0)
y_1(x),\label{cond1}\\
&&\varphi_2(x)=y'_1(0) y_2(x)-y'_2(0)
y_1(x),\label{cond2}\\
&&\varphi_{3}(x)=y_{1}''(0)y_{2}(x)-y_{2}''(0)y_{1}(x),\label{eq1}
\eean for $x\in \[-1,0\]$, and \bean
&&\psi_1(x)=z_1(0)z_2(x)-z_2(0)z_1(x),\label{cond3}\\
&&\psi_2(x)=z'_1(0)z_2(x)-z'_2(0)z_1(x),\label{cond4}\\
&&\psi_{3}(x)=z_{1}''(0)z_{2}(x)-z_{2}''(0)z_{1}(x),\label{eq2}
\eean for $x\in \[0,1\]$. Let $\phi_1(x,\la)$ and $\phi_2(x,\la)$ be
two solutions of Equation \eqref{s1} and \eqref{s2} which satisfy
the initial conditions $u(-1) = u''(-1)=0$ and $v(1) = v''(1)=0$,
respectively.
For clarity, the rest of the proof is divided into three steps.\\
\indent{\textsc{ STEP} $1$.} Assume that the subwronskians
 $\overline{\sigma}_i(0,\la)\neq 0$ for $i=1,2$.\\
 This implies that $\varphi'_1(0)\neq0$ and $\varphi_{2}(0)\neq0$
 (resp. $\psi'_1(0)\neq0$ and $\psi_{2}(0)\neq0$).
 Under this assumption $\varphi_i$ (resp. $\psi_i$), $i=1,2$, are linearly independent
 solutions of Equation \eqref{s1}(resp. Equation \eqref{s2}). Therefore, there exist constants $a$,
   $b$, $c$ and $d$ such that
   $$ \phi_1(x,\la)=a\varphi_1(x,\la)+b\varphi_2(x,\la)
   \hbox{ and }
   \phi_2(x,\la)=c\psi_1(x,\la)+d\psi_2(x,\la).$$
From the first two conditions of \eqref{s4}, we have \be
c=\frac{\varphi'_1(0)}{\psi'_1(0)}a\hbox{ and }
d=\frac{\varphi_2(0)}{\psi_2(0)}b. \label{alt3}\ee Using this
together with the last condition of \eqref{s4}, we get
   \be a\(\se_1\varphi''_1(0)-\frac{\se_2\varphi'_1\psi''_1(0)}{\psi'_1(0)}\)+b
   \(\se_1\varphi''_2(0)-\frac{\se_2\varphi_2\psi''_2(0)}{\psi_2(0)}\)=0.\ee
It is obvious that if
$\(\se_1\psi'_1\varphi''_1(0)-\se_2\varphi'_1\psi''_1(0)\)\neq0$ or
$\(\se_1\psi_2\varphi''_2(0)-\se_2\varphi_2\psi''_2(0)\)\neq0$, then
$\D E_{\la}=1$. Assume now the alternative case, i.e.,
\be\se_1\psi'_1\varphi''_1(0)-\se_2\varphi'_1\psi''_1(0)=0 \hbox{
and } \se_1\psi_2\varphi''_2(0)-\se_2\varphi_2\psi''_2(0)=0,
\label{al1}\ee then $a,~b\in\R.$ From \eqref{s5} and \eqref{alt3},
one gets \be\label{b1} a\(\mathcal{T}^1\varphi_1(0)-
\frac{\varphi'_1\mathcal{T}^2\psi_1(0)}{\psi'_1(0)}\)+b
\(\mathcal{T}^1\varphi_2(0)-
\frac{\varphi_2\mathcal{T}^2\psi_2(0)}{\psi_2(0)}\)=-M\la
b\varphi_2(0).\ee It can be easily verified from \eqref{eq:wrons1},
\eqref{eq:wronss1}, \eqref{cond1}-\eqref{cond2} and
\eqref{cond3}-\eqref{cond4}, that \be
\se_1\varphi''_2(0)=\overline\tau_1(0)=\mathcal{T}^1\varphi_1(0)
\hbox{ and } \se_2\psi''_2(0)=\overline
\tau_2(0)=\mathcal{T}^2\psi_1(0). \label{alt2}\ee Using this and
\eqref{al1}, we get
\begin{equation}\label{alt}
\(\mathcal{T}^1\varphi_1(0)-\frac{\varphi'_1\mathcal{T}^2\psi_{1}(0)}
{\psi'_1(0)}\)=\(\se_1\varphi''_2(0)-
\frac{\se_2\varphi_2\psi''_2(0)}{\psi_2(0)}\)=0.
\end{equation}
Now, if $b\neq0$, then by \eqref{alt}, Equality \eqref{b1} takes the
form \be\label{alt1}
\frac{\mathcal{T}^1\varphi_2(0)}{\varphi_2(0)}-\frac{\mathcal{T}^2\psi_2(0)}{\psi_2(0)}=-M\la.\ee
On the other hand, in view of Lemma \ref{Ls1m}, $\varphi_2$
satisfies one
of the following properties:\\
{\textit{Case} $1.$} $\varphi_2(0)>0$, $\varphi''_2(0)\geq0$ and
$\mathcal{T}^1
\varphi_2(0)\geq0$.\\
Then the first inequality in \eqref{alt} together with Lemma
\ref{Ls1m}, imply that $\psi_2\psi''_2(0)\geq0,$ and
$\psi_2\mathcal{T}^2\psi_2(0)<0$. As consequence, the left hand in
\eqref{alt1} is nonnegative, a
contradiction.\\
{\textit{Case} $2.$} $\varphi_2(0)<0$, $\varphi''_2(0)\leq0$ and
$\mathcal{T}^1
\varphi_2(0)\leq0$.\\
The {\textit{Case} $2$} can be handled in a same way.\\
{\textit{Case} $3.$} $\varphi_2\varphi''_2(0)<0$.\\
Then from \eqref{eq:wrons1} and \eqref{cond1}-\eqref{cond2} together
with \eqref{alt2}, we obtain
$$0>\se_1\varphi_2\varphi''_{2}(0)=-
\overline{\sigma}_1\overline\tau_1(0)=-\varphi_{1}'\mathcal{T}^1\varphi_{1}(0),$$
whence, $\varphi_{1}'\mathcal{T}^1\varphi_{1}(0)>0$. Again by Lemma
\ref{Ls1m}, we have $\varphi_{1}'\varphi_{1}''(0)>0$. According to
\eqref{al1} and \eqref{alt}, $\psi_{1}'\mathcal{T}^2\psi_{1}(0)>0$
and $\psi_{1}'\psi_{1}''(0)>0$. This is in
contradiction with Lemma \ref{Ls1m} and the condition \eqref{s3}.\\
Therefore $b=0$, and hence, from \eqref{alt3} we deduce that $\D
E_{\la}=1$.\\
\indent{\textsc{STEP} $2$.} Assume that $\overline\sigma_i(0)=0$ for $i=1,2$.\\
Under this assumption together with Lemma \ref{le:dim}, we have
$\varphi_1 '(0)=\psi_1'(0)=0$, $\varphi''_1(0)\neq0$ and
$\psi''_1(0)\neq0$. On the other hand by Lemma \ref{lem},
$\overline\sigma_i'(0)\neq0$ for $i=1,2$. This means that all of the
functions $\varphi_3$, $\varphi'_3$, $\psi_3$ and $\psi'_3$ does not
vanish at $x=0$. Therefore, $\varphi_1$ and $\varphi_3$ (resp.
$\psi_1$ and $\psi_3$) are linearly independent solutions of
\eqref{s1} (resp. of \eqref{s2}). Consequently, there exist
constants $a$,
   $b$, $c$ and $d$ such that
   $$ \phi_1(x,\la)=a\varphi_1(x,\la)+b\varphi_3(x,\la)
   \hbox{ and }
   \phi_2(x,\la)=c\psi_1(x,\la)+d\psi_3(x,\la).$$ Substituting these expressions
into the condition \eqref{s4}, one gets \be
b\frac{\varphi_3(0)}{\psi_3(0)}=d,~b\frac{\varphi'_3(0)}{\psi'_3(0)}=d~\hbox{
and }~c=\frac{\se_1\varphi_{1}''(0)}{\se_2\psi_{1}''(0)}a.
\label{alt8}\ee This implies that
$$b\(\frac{\varphi_3(0)}{\psi_3(0)}-\frac{\varphi'_3(0)}{\psi'_3(0)}\)=0.$$
Clearly, if $\(\psi_3\varphi'_3(0)-\varphi_3\psi'_3(0)\)\neq0$, then
$b=0$, and hence, $\D E_{\la}=1$. Now, suppose that
\begin{equation}\label{la3}
\psi_3\varphi'_3(0)=\varphi_3\psi'_3(0),
\end{equation}
then $a,~b\in \R$. A combination of \eqref{eq:wrons1},
\eqref{eq:wronss1}, \eqref{eq:wrons2} and \eqref{cond1}-\eqref{eq2}
yields
$$
\varphi_1''(0)=\overline{\si}_1'(0)=-\varphi_3(0),~~\mathcal{T}^1\varphi_1(0)=\overline{\tau}_1(0)=-\se_1\varphi'_3(0),
$$  and $$
\psi_1''(0)=\overline{\si}'_2(0)=-\psi_3(0),~~\mathcal{T}^2\psi_1(0)=\overline{\tau}_2(0)=-\se_2\psi'_3(0).
$$ Using these relations together with \eqref{la3}, one has
\be
\si_2\psi''_1\mathcal{T}^1\varphi_1(0)-\si_1\varphi''_1\mathcal{T}^2\psi_1(0)=0.\label{alt7}\ee
Since $\varphi_1,$ $\varphi_1'$ , $\psi_{1}$, and $\psi_{1}'$ vanish
at $x=0$, then by \eqref{s4} and \eqref{alt7}, the function
\begin{equation*}
\phi_0(x)=\left\{
              \begin{array}{ll}
                \si_2(0)\psi_{1}''(0)\varphi_1(x),&x\in\[-1,0\],\\
                \si_1(0)\varphi_{1}''(0)\psi_{1}(x),&x\in\[0,1\],
              \end{array}
            \right.
\end{equation*}
is a solution of the regular problem \eqref{s1}-\eqref{s5} for
$M=0$, which satisfies $\phi_0(0)=\phi_0'(0)=0$, and this is in
contradiction with Lemma \ref{le:cont}. Therefore, $b=0$, and hence
by \eqref{alt8}, we deduce that $\D E_{\la}=1$. \\
\indent{\textsc{STEP} $3$.} Assume that $\overline{\sigma}_1(0)=0$
and $\overline{\sigma}_2(0)\neq0$ (or conversely).\\ Then
$\varphi_1(0)=\varphi_{1}'(0)=0$. Let us recall from Step $1$, that
$\psi_1$ and $\psi_2$ are linearly independent solutions of
\eqref{s2}. Furthermore, from Step $2$,  $\varphi_1$ and $\varphi_3$
are linearly independent solutions of \eqref{s1}. Hence there exist
constants $a$,
   $b$, $c$ and $d$ such that
   $$ \phi_1(x,\la)=a\varphi_1(x,\la)+b\varphi_3(x,\la)
   \hbox{ and }
   \phi_2(x,\la)=c\psi_1(x,\la)+d\psi_2(x,\la).$$ Substituting these expressions
into the condition \eqref{s4}, we find \be
d=\frac{\varphi_3(0)}{\psi_2(0)}b,~c=\frac{\varphi'_{3}(0)}{\psi'_{1}(0)}b,\label{alt10}\ee
and \be
a\si_1\varphi''_{1}(0)=b\(\frac{\si_2\varphi'_{3}\psi''_{1}(0)}{\psi'_{1}(0)}+\frac{\si_2\varphi_3\psi''_{2}(0)}{\psi_2(0)}\).\label{alt9}\ee
Obviously, $\varphi_{1}''(0)\neq0$ since otherwise
$\overline{\sigma}_1'(0)=0$, and this is in contradiction with Lemma
\ref{lem}. Hence from \eqref{alt10} and \eqref{alt9}, $E_{\lambda}$
is generated by an unique solution, i.e., $\D E_{\lambda}=1$.\\
From the above, we deduce that the geometric multiplicity of the
eigenvalue $\la$ is equal to one. On the other hand, by Theorem
\ref{rr}, the linear operator $\mathcal{A}$ is self-adjoint in
$\mathcal{H}$, then $\la$ is algebraically simple. The proof is
complete.
\end{proof}\\
\begin{proof} {\bf of Assertion (b).}
It should be noted that by \eqref{s5}, $\mathcal{T}^1{u}(0)$ and
$\mathcal{T}^2{v}(0)$ do not necessarily coincide. Then it is not
possible to apply directly Lemma \ref{Ls1m}. Let
$\lambda_n,~(n\geq1)$ be an eigenvalue of Problem
\eqref{s1}-\eqref{s5} and $\phi_n$ be the corresponding
eigenfunction. Obviously,  $\phi_n$ can be written, up to a
multiplicative constant, in the unique form
$$ \phi_n(x,\la_n)= \left\{
\begin{array}{lll}
{u}_n(x)={u}(x,\la_n),&x\in\[-1,0\],\\
{v}_n(x)={v}(x,\la_n),&x\in\[0,1\].
\end{array}
\right.
$$
Suppose that
$${\phi}_{n}'\mathcal{T}^1{\phi}_{n}(-1)={u}_{n}'\mathcal{T}^1{u}_{n}(-1)\geq0 ~~\hbox{for some }n\geq1,$$
say ${u}_{n}'\geq0$ and $\mathcal{T}^1{u}_{n}(-1)\geq0$. Since
$u_n(-1)=u_n''(-1)=0$, then by the first statement of Lemma
\ref{Ls1m}, $u_n, u'_n, u''_n$ and $\mathcal{T}^1u_n$ are positive
at $x=0$. Thus from \eqref{s4}, one has
$$v_n(0)>0,~v'_n(0)>0\hbox{~~and~~} v''_n(0)>0.$$
Since $\la_n>0$, it follow from \eqref{s5}, that
$\mathcal{T}^2{v}_{n}(0)>0$. Again the first statement of Lemma
\ref{Ls1m} yields a contradiction with the boundary condition
$v_n(1)=0$. The proof of the second inequality in \eqref{alt12} is
similar.
\end{proof}
\addcontentsline{toc}{Section}{Bibliography}

\end{document}